\documentclass[]{preprint}
\usepackage{amsmath}
\usepackage{times}

\usepackage{mhequ}
\usepackage{mhenvs}
\usepackage{mhfig}
\usepackage{Symbols}

\makeatletter
\newtheorem{assumption}{Assumption}
\def\theassumption{H.\@arabic\c@assumption}
\makeatother

\newcommand{\T}{{\mathbf{T}^d}}
\newcommand{\e}{\varepsilon}
\renewcommand{\P}{\mathbf{P}}
\newcommand{\TV}{\mathrm{TV}}
\def\E{\mathbf{E}}
\def\cD{\mathcal{D}}
\def\bref#1{\textbf{\ref{#1}}}

\begin{document}

\title{Homogenization of periodic linear degenerate PDEs}
\author{Martin Hairer\inst{1}, Etienne Pardoux\inst{2}}
\institute{Department of Mathematics, The University of Warwick, Coventry CV4 7AL, United
Kingdom.
\email{hairer@maths.warwick.co.uk} \\
\and Laboratoire d'Analyse, Topologie, Probabilit\'es, 
Universit\'e de Provence, F-13453 Marseille, France.
\email{pardoux@cmi.univ-mrs.fr}}
\maketitle

\begin{abstract}\parindent1em
It is well-known under the name of `periodic homogenization' that, under a centering
condition of the drift, a periodic diffusion process on $\R^d$ converges, under diffusive rescaling, 
to a $d$-dimensional Brownian motion. Existing proofs of this result all rely on uniform
ellipticity or hypoellipticity assumptions on the diffusion. In this paper, we considerably weaken these
assumptions in order to allow for the diffusion coefficient to even vanish on an open set. 

As a consequence, it is no longer the case that the effective diffusivity matrix is necessarily non-degenerate. 
It turns out that, provided that some very weak regularity conditions are met, the range of the effective
diffusivity matrix can be read off the shape of the support of the invariant measure for the 
periodic diffusion. In particular, this gives some easily verifiable conditions for the 
effective diffusivity matrix to be of full rank. We also discuss the application of our
results to the homogenization of a class of elliptic and parabolic PDEs.
\end{abstract}

\section{Introduction}
Our goal is to study, by a probabilistic method, the limit as
$\eps\to0$ of the solution  $u^\eps(t,x)$ of an elliptic PDE in
the regular bounded domain $D\subset\R^d$
\begin{equation}\label{ellip}
\left\{
\begin{aligned}
{\mathcal L}_\eps u^\eps(x)+f\left(x,\frac{x}{\eps}\right)u^\eps(x)&=0,\, x\in D,\\
u^\eps(x)&=g(x),\, x\in\partial D,
\end{aligned}
\right.
\end{equation}
where $f$ is bounded from above, and $g$ is continuous, as well
as the limit of $u^\eps(t,x)$, the solution of a parabolic PDE of
the form
\begin{equation}\label{parab}
\left\{
\begin{aligned}
\frac{\partial u^{\eps}(t, x)}{\partial t}&={\mathcal L}_{\eps}u^{\eps}(t,x)
+\left(\frac{1}{\eps}e(\frac{x}{\eps})+f(x,\frac{x}{\eps})\right)
u^{\eps}(t,x)\;,
\\  u^{\eps}(0, x)&=g(x), x\in\R^d\;.
\end{aligned}  \right.
\end{equation}
In both cases, the linear operator $\CL_\eps$ is assumed to be a
second order differential operator with rapidly oscillating coefficients given by
\begin{equation*}
 {\mathcal L}_{\eps}=\frac{1}{2}\sum_{i,
j=1}^{d}a_{ij}\bigl(\frac{x}{\eps}\bigr)\frac{\partial^{2}}{\partial
x_{i}\partial x_{j}}
+\sum_{i=1}^{d}\Bigl[\frac{1}{\eps}b_{i}\bigl(\frac{x}{\eps}\bigr)
+c_i\bigl(\frac{x}{\eps}\bigr)\Bigr]
\frac{\partial}{\partial x_{i}}.
\end{equation*}
The novelty of our result lies in the fact that we allow the matrix $a$ to degenerate
(and even possibly to vanish)
in some open subset $D$ of $\R^d$. There is by now quite a vast literature
concerning the homogenization of second order elliptic and parabolic PDEs with a 
possibly degenerating matrix of second order coefficients $a$, 
see among others \cite{asc}, \cite{bebo}, \cite{bmt}, \cite{eppw}, \cite{fp}.
But, as far as we know, in all of these works, either the 
coefficient $a$ is allowed to degenerate in certain directions only, 
or else it may vanish
on sets of Lebesgue measure zero only. It seems that our paper presents the first results 
where the matrix $a$ is
allowed to vanish on an open set. The main technical difficulty that we have to overcome
is the lack of regularisation since we do not assume $\CL_\eps$ to be hypoelliptic (not even
on a set of full measure). However, it turns out that it is possible to show nevertheless that
under very weak assumptions, its resolvent maps $\CC^1$ into $\CC^1$ (see \lem{le2}),
which provides a $C^1$ solution to certain Poisson equations, and is sufficient to make an approximation argument work (see \lem{le:ito}).

Because of the high degree of degeneracy allowed by our approach, it is no longer
obvious that the effective diffusivity $A$ of the homogenized operator
\begin{equ}
\CL_0 = {1\over 2} \sum_{i,
j=1}^{d}A_{ij}\frac{\partial^{2}}{\partial
x_{i}\partial x_{j}}
+\sum_{i=1}^{d}C_i
\frac{\partial}{\partial x_{i}}
\end{equ}
is non-degenerate.
We shall therefore also seek to characterize the image of the homogenized diffusion matrix.
It turns out that this can be done in terms of the support of the invariant measure of the 
diffusion process on the torus $\T$ with drift $b$ and diffusion matrix $a$.

The paper is organized as follows. Section 2 contains our assumptions 
and several important preliminary results. Section 3 presents the homogenization result, in probabilistic terms. Section 4 contains our characterization of the image of the homogenized diffusion matrix $A$, and sections 5 and 6 present the application to elliptic and parabolic PDEs. Finally, Section~7 contains
a few concrete example that illustrate the scope of the results in this paper and highlight 
the differences with the existing literature.

\section{Assumptions and preliminary results}
Given
$\eps\geq 0$, $x\in \R^d$, let $\{X^{x,\eps}_t\}$ denote the solution of the
SDE
\begin{equation}\label{sde}
  X^{x,\eps}_{t}=x+
\int_{0}^{t}\Bigl[\frac{1}{\eps}b\left(\frac{X^{x,\eps}_{s}}{\eps}
\right) +c\left(\frac{X^{x,\eps}_{s}}{\eps}
\right)\Bigr]\,ds
+\sum_{j=1}^{m}\int_{0}^{t}\sigma_{j}\left(\frac{X^{x,\eps}_{s}}{\eps}
\right)\,dW^{j}_{s},
\end{equation}
where $b$, $c$, and $\sigma$  are periodic, of period one in each direction,
and the process $\{W_t=(W^1_t,\ldots,W^m_t),\ t\ge0\}$ is a standard $m$-dimensional Brownian motion. 

Define $\tilde X^{x,\eps}_{t}=\frac{1}{\eps}X^{x,\eps}_{\eps^{2}t}$.
Then there exists a standard $m$-dimensional Brownian motion
$\{W_t\}$, depending on $\eps$ (but we forget that dependence since it has no incidence on the
law of the process),
such that
\begin{equation}\label{sde_torus_eps}
\tilde X^{x,\eps}_{t}=\frac{x}{\eps}+\int_{0}^{t}\left[b(\tilde X^{x,\eps}_{s})+\eps c(\tilde X^{x,\eps}_{s})\right]\,ds
 +\sum_{j=1}^{m}\int_{0}^{t}\sigma_{j}(\tilde X^{x,\eps}_{s}
)\,dW^{j}_{s}
\end{equation}
In the sequel, we shall consider the solution of
\eqref{sde_torus_eps}, as taking values in the torus $\T$. We will also consider
the same equation starting from $x$, but without the term $\eps c$ in the drift, namely
 \begin{equation}\label{sde_torus}
\tilde X^{x}_{t}=x+\int_{0}^{t}b(\tilde X^{x}_{s})\,ds
 +\sum_{j=1}^{m}\int_{0}^{t}\sigma_{j}(\tilde X^{x}_{s})\,dW^{j}_{s}.
\end{equation}

We denote by $J_t^x$ the Jacobian of the stochastic flow associated to $\tilde X_t^x$, that
is $\{J^x_t,\ t\ge0\}$, the $d\times d$-matrix valued stochastic process solving
\begin{equ}[e:linear]
d J^x_t = Db(\tilde X_t^x)J_t^x\,dt + \sum_{j=1}^m D\sigma_j(\tilde X^x_t)J^x_t\,dW^j_t \;,
\quad J^x_0 = I\;.
\end{equ}

To the SDE satisfied by the process $\{\tilde{X}^x_{\cdot}\}$, we associate, inspired by 
Stroock-Varadhan's support theorem, the following controlled ODE ({\it from now on, we adopt the convention of summation over repeated indices}). For each $x\in\T$, 
$u\in L^2_{\text{loc}}(\R_+;\R^m)$, let $\{z_u^x(t),\ t\ge0\}$ denote the solution of
\begin{equation}\label{control}
\left\{
\begin{aligned}
\frac{dz_i}{dt}(t)&=b_i(z(t))
-\frac{1}{2}\left[\frac{\partial\sigma_{i j}}{\partial x_k}\sigma_{kj}\right](z(t))
+\sigma_{i j}(z(t))u_j(t),\\
z(0)&=x.
\end{aligned}
\right.
\end{equation}
We shall also need the controlled ODE with $b$ replaced by $b+\eps c$,
namely we shall denote by $\{z_u^{x,\eps}(t),\ t\ge0\}$
the solution of
\begin{equation}\label{control_eps}
\left\{
\begin{aligned}
\frac{dz_i}{dt}(t)&=\left[b_i+\eps c_i\right](z(t))
-\frac{1}{2}\left[\frac{\partial\sigma_{i j}}{\partial x_k}\sigma_{kj}\right](z(t))
+\sigma_{i j}(z(t))u_j(t),\\
z(0)&=x.
\end{aligned}
\right.
\end{equation}

We will throughout this paper make the following assumptions on the drift and the
diffusion coefficient
\begin{assumption}\label{h1}
The functions $\sigma$, $b$, and $c$ are of class $\CC^\infty$ and periodic of period $1$ in each
direction.
\end{assumption}

Consider now $b$ and the $\sigma_j$ as vector fields on the torus $\T$. We say that the \textit{strong H\"ormander
condition} holds at some point $x \in \T$ if the Lie algebra generated by $\{\sigma_j\}_{j=1}^m$
spans the whole tangent space of $\T$ at $x$. We furthermore say that the 
\textit{parabolic H\"ormander condition} holds at $x$ if the Lie algebra generated
by $\{\d_t + b\} \cup \{\sigma_j\}_{j=1}^m$ spans the whole tangent space of $\R \times \T$ at $(0,x)$.

\begin{assumption}\label{h2}
There exists a non-empty, open and connected subset $U$ of $\T$ on which the strong H\"ormander condition
holds. Furthermore, there exists $t_0 > 0$ and $\eps_0$ such that, for all $x \in \T$,
$0\le\eps\le\eps_0$, one has
\begin{equ}
\inf_{u\in L^2(0,t_0;\R^m)}\{\| u\|_{L^2};\ z_u^{x,\eps}(t_0)\in U\}<\infty\;.
\end{equ}
\end{assumption}
Note that, by upper semicontinuity, the supremum over $x\in\T$ and $0\le\eps\le \eps_0$ of the above 
infimum is bounded by a universal constant $K$.

Whenever $X$ is a random variable and $A$ an event, we shall use the notation $$\E(X;A) = \E(X\chi_A) = \E(X\,|\,A)\P(A).$$

\begin{assumption}\label{h4}
One has 
\begin{equ}
\inf_{t > 0} \sup_{x \in \T} \E \bigl(|J_t^x|\,;\,\{\tau_V^x\ge t\}\bigr)<1\;,
\end{equ}
where $V$ denotes the subset of $\T$ where the parabolic H\"ormander condition holds and 
$\tau_V^x$ is the first hitting time of $V$ by the process $\{\tilde{X}^x_t\}$.
\end{assumption}


\begin{remark}
One simple criteria for this assumption to hold is the existence of a time $t_2$ such that 
 $\P(\tau_V^x\le t_2)=1$ for all $x \in \T$.
\end{remark}

Denote by $p_\eps(t;x,A)$ the transition probabilities of the
$\T$-valued Markov process $\{\tilde X_t^{\eps,x}\}$. We shall write 
$p(t;x,A)$ for $p_0(t;x,A)$.
\begin{lemma}\label{le1}
Under the assumptions \bref{h1} and \bref{h2}, 
the following Doeblin condition is satisfied~: there exists
$t_1>0$, $0<\eps_1\le\eps_0$, $\beta>0$ and $\nu$ a probability measure on $\T$ which is absolutely continuous
with respect to Lebesgue's measure, s.t. for
all $0\le\eps\le\eps_1$, $x\in\T$, $A$ Borel subset of $\T$,
\begin{equation}\label{doeb}
p_\eps(t_1;x,A)\ge\beta\nu(A).
\end{equation}
\end{lemma}

\begin{proof}
It is well-known since the works of Malliavin, Bismut, Stroock et al \cite{Mal76SCV,Bismut81,KSAMI,KSAMII,Nor86SMC} 
that \bref{h1} and \bref{h2} imply
that the transition probabilities starting from $x \in U$ have a $\CC^\infty$ density
with respect to the Lebesgue measure, and that this density is also smooth in the initial
condition and in $\eps$. In particular this implies that, for every $t > 0$,  
$(\eps,x) \mapsto p_\eps(t;x,\cdot\,)$ is continuous from 
$[0,1]\times U$ into the set of probability measures on $\T$, equipped with the total variation distance. 
Let now $x_0 \in U$ be arbitrary. It follows that there exists an open neighbourhood
$U_0 \subset U$ of $x$, $0<\eps_1\le\eps_0$ and a probability measure $\nu$ such that 
$p_\eps(t;x,\cdot\,)\ge\frac{1}{2}\nu$ for every $0\le\eps\le\eps_1$, $x \in U_0$.
Since the measures $p_\eps(t;x,\cdot\,)$ are absolutely continuous, $\nu$ must also be
a.c. with respect to Lebesgue's measure.

Since $U_0 \subset U$ and since Assumption~\bref{h2} implies that \eref{control} is locally controllable
in $U$, it follows from the support theorem
that, for every $s > t_0$ and every $x\in\T$, $0\le\eps\le\eps_1$, one has $p_\eps(s;x,U_0) > 0$.
Since the Markov semigroup generated by the solutions of \eref{sde_torus} is Feller and
since $U_0$ is open, the function $(\eps,x) \mapsto p_\eps(s;x,U_0)$ is lower semicontinuous
and therefore attains its lower bound $\beta'$ on $[0,\eps_1]\times\T$. The claim follows by taking
$\beta = \beta'/2$ and $t_1 = s+t$.
\end{proof}

In the sequel we shall denote by $\{\CP_{\e,t},\ t\ge 0\}$ the Markov semigroup associated to the $\T$-valued diffusion process 
$\{\tilde X^{\e,x}_t,\ t\ge 0\}$. We shall write $\CP_t$ for $\CP_{0,t}$.


\begin{corollary}\label{cor01}
Under the assumptions \bref{h1} and \bref{h2},  for each $0\le\eps\le\eps_1$, $\{\tilde
X_t^\eps\}$ possesses a unique invariant probability measure $\mu_\eps$, which is absolutely
continuous with respect to Lebesgue's measure. Moreover, there exist constants $C$ and
$\rho>0$ such that
\begin{equ}[0_spectral_gap]
\sup_{x\in\T,\ 0\le\eps\le\eps_1}\|p_\eps(t;x,\cdot)-\mu_\eps\|_{TV}\le C e^{-\rho t}\;,
\end{equ}
for every $t \ge 0$.
\end{corollary}

\begin{proof}
It follows from \eref{doeb} that for all $0\le\eps\le\eps_1$
\begin{equ}[contract]
\|\CP_{\eps,t_1} \mu_1 - \CP_{\eps,t_1} \mu_2\|_\TV \le (1-\beta) \|\mu_1 - \mu_2\|_\TV
\end{equ}
for any two probability measure $\mu_1$ and $\mu_2$ on $\T$.
This immediately implies the
existence and uniqueness of an invariant probability measure $\mu_\eps$ for $\tilde X_t^\eps$.

To see that $\mu_\eps$ is absolutely continuous with respect to Lebesgue's measure, note that one
can decompose the transition semigroup $\CP_{\eps,t}$ for $t = t_1$ as
\begin{equ}[e:exprmu]
\CP_{\eps,t_1} = \beta \nu + (1-\beta) \bar \CP_\eps
\end{equ}
for some Markov operator
 $\bar \CP_\eps$. 
It follows from the fact that the solutions of \eref{sde_torus_eps}
generate flows of diffeomorphisms \cite{Kun90}, that
 $\bar \CP^n_\eps \nu$ is absolutely continuous with respect to the Lebesgue measure for
every $n$. 
Fix now an arbitrary
set $A$ with Lebesgue measure $0$. It follows from the invariance of 
$\mu_\eps$ that
$\mu_\eps = \bigl(\CP_{\eps,t_1}^*\bigr)^n \mu_\eps$
for every $n > 0$. From this, \eref{e:exprmu}, and the absolute continuity of $\nu$, it follows 
that $\mu_\eps(A) \le (1-\beta)^n$ for every $n$,
and therefore that $\mu_\eps(A) = 0$.

Finally \eref{0_spectral_gap} follows from iterating \eref{contract} with $\mu_2=\mu_\eps$.
\end{proof}
We shall need the
\begin{lemma}\label{le:conv-mesinv}
Denote $\mu_0$ by $\mu$. As $\eps\to0$,
$$\mu_\eps\to\mu\quad\text{weakly}.$$
\end{lemma}
\begin{proof}
The tightness is obvious, since $\T$ is compact. Hence from any sequence
$\eps_n$ converging to 0, we can extract a further subsequence, which we still denote by $\{\eps_n\}$, such that 
$$\mu_{\eps_n}\to\tilde\mu.$$
Now if $f\in\CC(\T)$, $t>0$, clearly $\CP_{\eps_n,t}f(x)\to\CP_tf(x)$
as $n\to\infty$, uniformly for $x\in\T$, hence
$$\int_\T \CP_{\eps_n,t}f(x)\mu_{\eps_n}(dx)\to\int_\T\CP_tf(x)\tilde\mu(dx).$$
But the left hand side equals
$$\int_\T f(x)\mu_{\eps_n}(dx)\to\int_\T f(x)\tilde\mu(dx).$$
Consequently $\tilde\mu$ is invariant under $\CP_t^\ast$ for all $t>0$, hence $\tilde\mu=\mu$, and $\mu_\eps\to\mu$ weakly, as $\eps\to0$.
\end{proof}
We can now deduce from \eref{0_spectral_gap}, using also Lemma
\ref{le:conv-mesinv}, exactly as in Proposition 2.4 of \cite{PARDjfa99},
the
\begin{corollary}\label{cor02}
Whenever $f\in L^\infty(\T)$, for any $t>0$,
$$\int_0^tf\left(\frac{X^{x,\eps}_s}{\eps}\right)\,ds
\to t\int_\T f(y)\mu(dy)$$
in probability, as $\eps\to0$.
\end{corollary}

We finally assume that
\begin{assumption}\label{h5}
The drift $b$ satisfies the centering condition
$\int_\T b(x)\mu(dx)=0$.
\end{assumption}

Denoting by 
$$L=\frac{1}{2}a_{ij}(x)\frac{\partial^{2}}{\partial x_{i}x_{j}}+b_{i}(x)
\frac{\partial}{\partial x_{i}}$$
with $a(x)=\sigma\sigma^\ast(x)$
the infinitesimal generator of the $\T$-valued diffusion process 
$\{\tilde{X}^x_t,\ t\ge0\}$, it follows from Lemma~\ref{le2} 
below that under assumptions \bref{h1}--\bref{h5},
there exists a unique 
$\CC^1(\T;\R^d)$ solution of the Poisson equation
$$L\widehat{b}(x)+b(x)=0\;,\quad x\in \T\;.$$
This solution is given by \cite{ep_av3}
\begin{equ}
\widehat b (x) = \int_0^\infty \CP_t b(x)\,dt\;.
\end{equ}

It is not clear \textit{a priori} that $\widehat b$ is differentiable, however the following result
shows that since $b$ belongs to $\CC^1(\T)$, so does $\widehat b$.
\begin{lemma}\label{le2}
Under assumption \bref{h1}, the semigroup $\CP_t$ generated by \eref{sde_torus}
maps $\CC^1(\T)$ into itself. If furthermore \bref{h2} and \bref{h4} hold, then there
exist positive constants $C$ and $\gamma$ such that
\begin{equ}[e:gap]
\|\CP_t f\|_{\CC^1(\T)} \le Ce^{-\gamma t}\| f\|_{\CC^1(\T)}\;,
\end{equ}
for every $f \in \CC^1(\T)$ such that $\int f(x)\mu(dx) = 0$,
and for every $t > 0$. 
\end{lemma}
\begin{proof}
The proof relies on the techniques developed in
\cite{HaiMat06SGW}.
Suppose that we can find a time $t_2$ such that there exists a constant $C$ and $\delta > 0$ with
\begin{equ}[e:semigap]
\|D \CP_{t_2} f\|_{\L^\infty} \le C \|f\|_{\L^\infty} + {(1-\delta)} \| Df\|_{\L^\infty}\;,
\end{equ}
for every test function $f \in \CC^1(\T)$.
It follows from \cite[Section~2]{HaiMat06SGW} that this, together with the Doeblin condition 
given by Lemma~\ref{le1}, implies that \eref{e:gap} holds.

It remains to show that \eref{e:semigap} does indeed hold. 
The set $V$ can be defined as
$$V=\Big\{x;\ Q(x):=\sum_{Y\in\mathcal{Y}} Y(x)\otimes Y(x)>0\Big\},$$
where $\CY$ denotes the set of vector fields consisting of the $\sigma_j$'s, their brackets and their brackets with $b$, the number of those
being bounded by $d$. The inequality $Q(x)>0$ means that the quadratic form $Q(x)$ is strictly positive definite.
We let $V^k$ be the
set of points where $Q(x)>k^{-1} I$.  

We now want to take advantage of the fact that as soon as the process 
$\tilde{X}^x$ hits $V_k$, its Malliavin matrix becomes invertible.
However, we face the problem that
the hitting time of $V_k$ doesn't necessarily have a Malliavin derivative. This is the motivation for the next construction. Let $\Omega$ and $\bar\Omega$ be two independent copies 
 of the $d$-dimensional Wiener 
space, and denote by $B$ and $\bar B$ the corresponding canonical processes. 
We first consider the auxiliary SDE
\begin{equ}
dy = b(y)\,dt + \sum_j \sigma_j(y)\,dB\;,\qquad y(0) = {x}\;,
\end{equ} 
and define $\tilde\tau^x_k$ as the hitting time of $V_k$ by the process 
$y$. Let $t_2>0$ be such that
\begin{equ}[choice:t2]
\sup_{x \in \T} \E \bigl(|J_{t_2}^x|\,;\,\{\tau_V^x\ge t_2\}\bigr)<1,
\end{equ}
and
$$\tau^x_k=\begin{cases} \tilde\tau^x_k,&\text{when $\tilde\tau^x_k< t_2$};\\
                      +\infty,&\text{otherwise}.\end{cases}$$
 We now define a
 process $W$ by
\begin{equ}
W_t = \left\{\begin{array}{cl} B_t & \text{for $t \le \tau^x_k$,} \\ B_{\tau^x_k} + \bar B_{t-{\tau^x_k}} & \text{for $t \ge {\tau^x_k}$.} \end{array}\right.
\end{equ}
Since ${\tau^x_k}$ is a stopping time, $W$ is again a $d$-dimensional Wiener process.
We can (and will between now and the end of this proof) therefore consider the process
$\tilde X^x_t$ to be driven by the process $W$ that was just constructed.

Consider (see \eref{e:linear}) the $d\times d$-matrix valued SDE
\begin{equ}
d J^x_t = Db(\tilde{X}^x_t)J^x_t\,dt + \sum_{j=1}^m 
\bigl(D\sigma_j(\tilde X^x_t)J^x_t\bigr)\,dW^j_t \;,
\quad J_0 = I\;.
\end{equ}

We next define the (random) linear map 
$A_x\colon \L^2([0,1];\R^m) \to \R^d$ by
\begin{equ}
A_x u = \sum_{j=1}^m\int_{\tau^x_k}^{{\tau^x_k} + 1} (J^x_s)^{-1} \sigma_j(X^x_s) u_j(s-\tau^x_k)\,ds\;,
\end{equ}
with the convention that $A_x = 0$ on the set $\{{\tau^x_k} = +\infty\}$.
Define the (random) map $C_x \colon \R^d \to \R^d$ by
$C_x = A_x A_x^*$.  It then follows from \cite{Nua95} and the strong Markov property that
\begin{proposition}\label{prop:boundv}
For every $x \in \T$, the matrix $C_x$ is almost surely invertible on the set $\{{\tau^x_k} < +\infty\}$.
 Furthermore, for every $k$ and every $p$ there exists a constant $K_{k,p}$ such that
$\E\bigl(\|C_x^{-p}\|\,|\, {\tau^x_k} < \infty\bigr) \le K_{k,p}$.
\end{proposition}
Fix now an arbitrary vector $\xi \in \R^d$ and define a stochastic process $v\colon [0,1]\to\R^m$ by
\begin{equ}
v = A_{x}^* C_{x}^{-1} \xi\;.
\end{equ}
(We set as before $v = 0$ on the set $\{{\tau^x_k} = +\infty\}$.)
The reason for performing this construction is that we are now going to use Malliavin calculus
on the probability space $\Sigma = \Omega \times \bar\Omega$ equipped with 
the canonical Gaussian subspace generated by the Wiener process $\bar B$.
It is then straightforward to check \cite{Nua95} that one has $v \in D^{1,2}$. Note however that since
the stopping time ${\tau^x_k}$ is not Malliavin differentiable in general, it is \textit{not} 
necessarily true that  $v\in D^{1,2}$ if we equip $\Sigma$ with the Gaussian structure inherited 
from $W$.

This construction yields the following equality, valid for every
$f\in\CC^1(\T)$ and for $t = t_2+1$:
\begin{equ}
\int_0^t \cD_s \bigl(f(\tilde X^x_t)\bigr)\, v_s\, ds = \left\{\begin{array}{cl} (Df)(\tilde X^x_t) J^x_t \xi & \text{if ${\tau^x_k} < \infty$,} \\ 0 & \text{if ${\tau^x_k} = \infty$,} \end{array}\right. 
\end{equ}
where $\cD$ stands for the Malliavin derivative with respect to $\bar B$.
It then follows from the integration by parts formula \cite[Sect.~1.3]{Nua95} that
\begin{equ}
D \CP_t f(x) \xi = \E \Bigl(f(\tilde X^x_t) \int_0^1 v(s)\,d\bar B(s)\Bigr) - \E\bigl( (Df)(\tilde X^x_t) J^x_t \xi \,;\, {\tau^x_k} = \infty\bigr)\;,
\end{equ}
where the integral is to be understood in the Skorokhod sense. It follows from \cite[p.~39]{Nua95} that
\begin{equs}
\| D\CP_{t_2} f\|_\infty &\le
\|Df\|_\infty \E\bigl( |J^x_{t_2}|;\{\tau^x_k = \infty\}\bigr)  \\
&\quad + \|f\|_\infty \E \Bigl(\int_0^1 \|v_s\|^2\,ds +
\int_0^1 \int_0^1 \|\cD_s v_r \|^2\, dr\, ds\Bigr)\;.
\end{equs}
Since, for every $k$, the second term is bounded uniformly in $x\in\T$, it remains to show that we can choose $k$ such that
$$\sup_x \E\bigl( |J^x_{t_2}|;\{\tau^x_k = \infty\}\bigr)<1.$$
Assume this is not the case. Then to each $k$, we can associate a point
$x_k$ such that
$$\E\bigl( |J^{x_k}_{t_2}|;\{\tau^{x_k}_k \ge t_2\}\bigr)\ge1.$$
We can assume (up to extracting a subsequence) that $x_k\to x$. For
all $\ell\in\N$,
\begin{equs} 
1 \le \limsup_{k\to\infty}\E\bigl(|J^{x_k}_{t_2}|;\{\tau^{x_k}_k \ge t_2\}\bigr)&\le
\limsup_{k\to\infty}\E\bigl(|J^{x_k}_{t_2}|;\{\tau^{x_k}_\ell \ge t_2\}\bigr)\\
&\le\E\bigl(|J^x_{t_2}|;\{\tau^x_\ell\ge t_2\}\bigl),
\end{equs}
since the mapping $x\to\chi_{\{\tau^x_\ell\ge t_2\}}$ is a.s. upper semicontinuous,
as the indicator function of a closed subset of trajectories.  Since
$$\lim_{\ell\to\infty}\chi_{\{\tau^x_\ell\ge t_2\}}=\chi_{\{\tau^x_V\ge t_2\}}\ \text{ a.s.},$$
this contradicts \eref{choice:t2}.
\end{proof}
\begin{remark}
We claim that the assumption \bref{h4} is close to being sharp for the 
conclusion of Lemma \ref{le2} to hold. The simplest example (which is however
not a diffusion on the torus) on which this can be seen is as follows.
Let $\tau$ be an exponential random variable with parameter $b>0$, and 
$$X^x_t=\begin{cases} e^{at}x,&\text{if $t<\tau$};\\
                       0,&\text{otherwise}.\end{cases}  $$
In this case
$$\E\bigl(|J^x_t|;\ \{\tau>t\}\bigl)=e^{(a-b)t}\;,$$
and
$$\|D\CP_tf\|_\infty=e^{(a-b)t}\|Df\|_\infty\;,$$
so that $\|\CP_tf\|_{\CC^1} \to 0$ if and only if \bref{h4} holds.
\end{remark}
\section{The homogenization result}
The goal of this section is to show the
\begin{theorem}\label{th:homo} 
We have that, in the sense of weak convergence on the space $\CC(\R_+)$
equipped with the topology of uniform convergence on compact sets,
\begin{equ}
X^{x,\eps}\Rightarrow X^x\;,\quad\text{where}\quad
X^x_t=x+Ct+A^{1/2}W_t\;,
\end{equ}
as $\eps \to 0$. Here,
$\{W_t,\ t\ge0\}$ is a standard Brownian motion, and the homogenized
coefficients $C$ and $A$ are given by
\begin{align*}
C&= \int_{\T}(I+\nabla\widehat{b})c(x)\ \mu(dx)\;, \\
A&=\int_{\T}(I+\nabla\widehat{b})a(I+\nabla\widehat{b})^{*}(x)\ \mu(dx)\;.
\end{align*}
\end{theorem}
In order to prove this theorem, we need to get rid of the term of order
$\eps^{-1}$ in the SDE \eref{sde}. The trick is to replace $X^{x,\eps}_t$ by  $$\hat X^{x,\eps}_t:=X^{x,\eps}_t+\eps\hat b\left(\frac{X^{x,\eps}_s}{\eps}\right),$$
and to take advantage of the
\begin{lemma}\label{le:ito}
The following equality holds almost surely:
$$\hat{X}^{x,\eps}_t =x+\eps\hat{b}\left(\frac{x}{\eps}\right)
+\int_0^t(I+\nabla\hat{b})c\left(\frac{X^{x,\eps}_s}{\eps}\right)\,ds
+
\int_0^t(I+\nabla\hat{b})\sigma\left(\frac{X^{x,\eps}_s}{\eps}\right)\,dW_s\;.
$$
\end{lemma}
\begin{proof}
Let $\rho~:\R^d\to\R_+$ be a smooth
function with compact support, such that
$$\int_{\R^d}\rho(x)\,dx=1,$$ and define $\rho_n(x)=n^d\rho(nx)$. We
 regularize $\hat{b}$ by convolution
$$\hat{b}_n=\hat{b}\ast\rho_n.$$
We now deduce from It\^o's formula that
\begin{equation*}
\begin{split}
\hat{X}^{\eps,n}_t&:=X^{x,\eps}_t+\eps\hat{b}_n\left(\frac{X^{x,\eps}_t}{\eps}\right) \\
&=x+\eps\hat{b}_n\left(\frac{x}{\eps}\right)
+\eps^{-1}\int_0^t(L\hat{b}_n+b)\left(\frac{X^{x,\eps}_s}{\eps}\right)\,ds
+
\int_0^t(I+\nabla\hat{b}_n)c\left(\frac{X^{x,\eps}_s}{\eps}\right)\,ds\\
&\quad +
\int_0^t(I+\nabla\hat{b}_n)\sigma\left(\frac{X^{x,\eps}_s}{\eps}\right)\,dW_s.
\end{split}
\end{equation*}
We want next to let $n\to\infty$. Clearly $\hat{b}_n\to\hat{b}$
and $\nabla\hat{b}_n\to\nabla\hat{b}$ pointwise, and the two
sequences are bounded, uniformly with respect to $n$ and $x\in\T$.
It remains to treat the term containing the second order derivatives. One has
$$L\hat{b}_n=(L\hat{b})\ast\rho_n+\varphi_n,$$
 and since $\hat{b}$ is a weak solution of the Poisson
 equation,
$$(L\hat{b})\ast\rho_n=-b\ast\rho_n\to-b,$$
the sequence being again uniformly bounded. It remains to study
the sequence $\varphi_n$. Using again the convention of summation over
repeated indices, we have
\begin{equation}\label{fi}
\begin{split}
\varphi_n(x)&=\frac{1}{2}\int_{\R^d}[a_{ij}(x)-a_{ij}(x-y)]\frac{\partial^2\hat{b}}{\partial
x_i\partial x_j}(x-y)n^d\rho(ny)dy\\&\ +
\int_{\R^d}[b_i(x)-b_i(x-y)]\frac{\partial\hat{b}}{\partial
x_i}(x-y)n^d\rho(ny)dy.
\end{split}
\end{equation}
The fact that the second integral in \eqref{fi} converges to 0 and
is uniformly bounded is easily established. We now integrate by
parts the first integral, yielding
\begin{equation*}
\begin{split}
&\int_{\R^d}[a_{ij}(x)-a_{ij}(x-y)]\frac{\partial^2\hat{b}}{\partial
x_i\partial x_j}(x-y)n^d\rho(ny)dy\\
&\,=\int_{\R^d}\frac{\partial a_{ij}}{\partial
x_i}(x-y)\frac{\partial\hat{b}}{\partial x_j}(x-y)n^d\rho(ny)dy\\
&\quad +\int_{\R^d}[a_{ij}(x)-a_{ij}(x-y)]\frac{\partial\hat{b}}{\partial
x_j}(x-y)n^{d+1}\rho'_i(ny)dy
\end{split}
\end{equation*}
The first term on the right hand side is uniformly bounded and
converges to
$$\frac{\partial a_{ij}}{\partial
x_i}(x)\frac{\partial\hat{b}}{\partial x_j}(x)\;.$$ The second
term is equal to
$$\int_{\R^d}y\cdot\nabla a_{ij}(x-y')
\frac{\partial\hat{b}}{\partial x_j}(x-y)n^{d+1}\rho'_i(ny)dy\;,
$$
where $|y'|\le|y|$. This last quantity is equal to a bounded
sequence converging to 0, plus
$$\frac{\partial a_{ij}}{\partial x_k}(x)\frac{\partial\hat{b}}{\partial
x_j}(x)\int y_kn^{d+1}\rho'_i(ny)dy=-\frac{\partial
a_{ij}}{\partial x_i}(x)\frac{\partial\hat{b}}{\partial x_j}(x).$$
The lemma is established.
\end{proof}
We now proceed with the
\begin{proof}[of Theorem \ref{th:homo}]
We first note that since 
$$\left|\hat X^{x,\eps}_t-X^{x,\eps}_t\right|\le C\eps, $$
the Theorem will follow if we prove that, as $\eps \to 0$,
$$ \hat X^{x,\eps}\Rightarrow X^x.$$  
But since
\begin{equ}
\hat X^{x,\eps}_t=x+\int_0^t(I+\nabla\hat b)c\left(
\frac{X^{x,\eps}_s}{\eps}\right)\,ds+M^\eps_t\;,
\end{equ}
one has
\begin{equ}
\langle\!\langle M^{\eps}\rangle\!\rangle_t =
\int_0^t(I+\nabla\hat b)a(I+\nabla\hat b)^\ast\left(
\frac{X^{x,\eps}_s}{\eps}\right)\,ds\;,
\end{equ}
and the result follows from Corollary~\ref{cor02} and the martingale central limit theorem, see e.g. 
\cite[Thm~7.1.4]{EthKur86MP}.
\end{proof}
We conclude this section with the following result, which extends the averaging
result of \cor{cor02}. It is going to be needed in the applications of sections 5 and 6.
\begin{proposition}\label{pr:conv-erg}
Let $f\in \CC(\R^d\times\T)$. Then for every $t>0$,
the following convergence holds in law as $\eps\to0$~:
$$\int_0^tf\left(X^{x,\eps}_s,\frac{X^{x,\eps}_s}{\eps}\right)\,ds
\quad\Rightarrow\quad \int_0^t\overline f(X^x_s)\,ds,$$
where
$$\overline f(x):=\int_\T f(x,y)\,\mu(dy).$$
\end{proposition}
\begin{proof}
It is easily checked that $\overline f$ is continuous, hence
$$ \int_0^t\overline f(X^{x,\e}_s)\,ds \quad\Rightarrow\quad  
\int_0^t\overline f(X^{x}_s)\,ds, $$
as $\e\to0$. It then suffices to show that
$$\left|\int_0^tf\left(X^{x,\eps}_s,\frac{X^{x,\eps}_s}{\eps}\right)\,ds
-\int_0^t\overline f(X^{x,\e}_s)\,ds\right|\to 0$$
in probability, as $\e\to0$. This is proved by
 exploiting the tightness in $\CC(\R_+)$ of 
the collection of processes $\{X^{x,\eps}_\cdot,\ \eps>0\}$, and
Corollary \ref{cor02}. For the details, see the proof of
Lemma 4.2 in \cite{PARDjfa99} (there is a misprint in the statement of that Lemma).
\end{proof}
\section[The image of the homogenized diffusion matrix]{Characterization of the image of the homogenized diffusion matrix}
The aim of this section is to give a characterisation of the range of $A$ which can
very easily be read off the support of the invariant measure $\mu$. We first state
all of our results and then prove them in the following subsection.

\subsection{Statements of the results}
It will be convenient in our statements to choose an arbitrary particular point
$x_0\in U$. We define the set $\Gamma$, consisting of loops on $\T$, starting and ending at $x_0$, which take the form
$$\gamma(s)=z_u^{x_0}(s),\quad 0\le s\le t,$$
with arbitrary $t>0$ and $u\in L^2(0,t;\R^m)$. We call loops in $\Gamma$ \textit{admissible}.
It follows from \bref{h2} that the set of admissible loops does not depend on the particular
choice of $x_0$.

We next define a mapping $g : \Gamma\to \Z^d$ as follows. Lifting each loop $\gamma$,
so as to transform it into a curve $\bar{\gamma}$ in $\R^d$, we define
$$g(\gamma)=\bar{\gamma}(t)-\bar\gamma(0),$$ 
where $t$ is the end time of the loop $\gamma$. The first step in our characterization
of the image of $A$ is to show that

\begin{theorem}\label{th1} Under the 
assumptions \bref{h1}--\bref{h5}, $g(\Gamma)\subset\text{Im}A$.
\end{theorem}

The converse to this result takes the following form:

\begin{theorem}\label{th2} Under the 
assumptions \bref{h1}--\bref{h5},
for each $e\in\text{Im}A$, with $|e|=1$, and for every $\delta>0$, there exists
a loop $\gamma\in\Gamma$ such that
$$\left|\frac{g(\gamma)}{|g(\gamma)|}-e\right|\le\delta.$$
\end{theorem}
It follows immediately from the above theorems that under the same assumptions,
\begin{corollary}\label{cor1} 
$\Im A=\mathrm{span} \{ g(\gamma),\ \gamma\in\Gamma\}$.
\end{corollary}
\begin{corollary}\label{cor2}
The matrix $A$ is non degenerate if and only if there exists a collection of admissible loops $\{\gamma_1,\ldots,\gamma_d\}
\subset\Gamma$ such that
$$\mathrm{span}\{g(\gamma_1),\ldots,g(\gamma_d)\}=\R^d.$$
\end{corollary}
The next corollary is slightly less obvious and will be shown in the next subsection. It will mainly
be useful in the proof of \theo{th3} below, but since it is not \textit{a priori} an obvious fact,
we state it separately.
\begin{corollary}\label{cor3}
The set $G=\{g(\gamma),\ \gamma\in\Gamma\}$ is a subgroup of $\Z^d$.
\end{corollary}
We denote by $S$ the interior of the support of the invariant measure $\mu$.
It follows from \bref{h2} that the support of $\mu$ is then the closure of $S$.
It turns out that as far as the characterisation of $\Im A$ is concerned, we can replace 
admissible loops by arbitrary loops in $S$. This may sound somewhat surprising at first sight,
since it is certainly not true in general that every loop in $S$ is admissible. One has however

\begin{theorem}\label{th3}
Every admissible loop $\gamma$ satisfies $\gamma\subset S$. Conversely,
for every loop $\rho\subset S$, one can find an admissible loop $\gamma$ with 
$g(\rho)=g(\gamma)$.
\end{theorem}

\begin{remark}
Combining \theo{th3} with \cor{cor1}, we find that $\Im A$ is completely characterised by the topology of
the set $S$. This is however not necessarily true if we drop the regularity conditions \bref{h2}--\bref{h4}.
In particular, it follows from \cor{cor1} that one can find a basis for $\Im A$ which has rational coordinates. A `counterexample' (with $d=2$ and $m=1$) to this claim is given by taking $b = c = 0$
and $\sigma_1 = (1, \sqrt2)$, so that $\Im A$ is the vector space generated by $(1,\sqrt 2)$. Of course,
this example satisfies neither \bref{h2} nor \bref{h4} and does therefore not contradict our results.
\end{remark}

\begin{remark}\label{rem:JKO}
The argument on page 264 in \cite{JKO} shows that $A$ is non degenerate
whenever the set span$\{g(\gamma),\ \gamma\in S_0\}$ is equal to $\R^d$, where $S_0$
is a connected subset of the support of the invariant measure $\mu$ where $a$ is elliptic. Clearly the result presented here is stronger.
\end{remark}

\subsection{Proofs of the results}

We start with the the
\begin{proof}[of Theorem \ref{th1}]
Let $\xi\in\R^d$ be such that $\langle A\xi,\xi\rangle=0$. We will prove that for any $\gamma\in\Gamma$, $\langle g(\gamma),\xi \rangle=0$. We shall make use of the notation $b_\xi(x)=\langle b(x),\xi\rangle$,
and of the three following facts~:
\begin{itemize}
\item
$\widehat{b}_\xi$ solves the Poisson equation
\begin{equation}\label{poisson}
L\widehat{b}_\xi(x)+b_\xi(x)=0,\quad x\in\T;
\end{equation}
\item $\langle A\xi,\xi\rangle=0$, which, when we explicit the matrix $A$, amounts to
the fact that
\begin{equation}\label{ident2}
\langle\nabla\widehat{b}_\xi,\sigma_{\cdot j}\rangle+\langle\sigma_{\cdot j},\xi\rangle=0,\
1\le j\le d,\ \text{for }\mu \text{ a.\ e. }x.
\end{equation}
Note that \eqref{ident2} holds true $\mu$-almost everywhere. It can be noted that the trajectories
of the solutions of \eqref{control}
starting from the point $x_0$ remain in the interior of the support of $\mu$, which
is absolutely continuous with respect to the Lebesgue measure. Hence \eqref{ident2}
holds true almost everywhere in the neighbourhood of such a trajectory, consequently everywhere there, due to continuity. 
\item
We can now differentiate \eqref{ident2} with respect to $x_k$, and multiply the
resulting identity by $\sigma_{kj}$, from which we deduce that
\begin{equation}\label{ident3}
\sigma_{kj}(x)\frac{\partial}{\partial x_k}\left(\sigma_{ij}\frac{\partial\widehat{b}_\xi}{\partial x_i}\right)(x)=-\xi_i\left(\sigma_{kj}\frac{\partial\sigma_{ij}}{\partial x_k}
\right)(x),\quad x\in\T.
\end{equation}
\end{itemize}

Choose a loop $\gamma\in\Gamma$, and denote by $t$ its end time. We now have, for $0<s<t$, 
\begin{align*}
\frac{d}{ds}\widehat{b}_\xi(\gamma(s))&=\nabla\widehat{b}_\xi(\gamma(s))\frac{d\gamma}{ds}(s)\\
&=\langle b,\nabla\widehat{b}_\xi\rangle(\gamma(s))
-\frac{1}{2}\left[\frac{\partial\widehat{b}_\xi}{\partial x_i}
\frac{\partial\sigma_{ij}}{\partial x_k}\sigma_{kj}\right]\!(\gamma(s))
+\left[\frac{\partial\widehat{b}_\xi}{\partial x_i}\sigma_{ij}\right]\!(\gamma(s))u_j(s)\\
&=-b_\xi(\gamma(s))
-\frac{1}{2}\left[a_{ij}\frac{\partial^2\widehat{b}_\xi}{\partial x_i\partial x_j}
+\sigma_{kj}\frac{\partial \sigma_{ij}}{\partial x_k}\frac{\partial\widehat{b}_\xi}{\partial x_i}\right]\!(\gamma(s))-\xi_i\sigma_{ij}(\gamma(s))u_j(s)\\
&=-b_\xi(\gamma(s))+\frac{1}{2}\xi_i\left[\sigma_{kj}\frac{\partial\sigma_{ij}}{\partial x_k}\right]\!(\gamma(s))-\xi_i\sigma_{ij}(\gamma(s))u_j(s)\\
&=-\scal[B]{\frac{d\gamma}{ds}(s),\xi},
\end{align*}
where we have used \eqref{poisson} and \eqref{ident2} for the third equality,
and \eqref{ident3} for the fourth equality.
Integrating from $s=0$ to $s=t$, we deduce that
$$\langle g(\gamma),\xi\rangle+\langle\widehat{b}_\xi(\bar\gamma(t))-\widehat{b}_\xi(\bar\gamma(0)),\xi\rangle=0,
$$
from which it follows, since $\widehat{b}$ is periodic, that $\langle g(\gamma),\xi\rangle=0$. 
The result follows.
\end{proof}

\begin{proof}[of Theorem \ref{th2}]
Let $e\in\text{Im}A$ with $|e|=1$ and $\alpha>0$. We denote by $B(x,\alpha)$
the ball in $\R^d$ of radius
$\alpha$, centered at $x$. Note that
$$\P(A^{1/2}W_1\in B(e,\alpha))>0.$$
It follows that for $\e>0$ small enough (here $X^{\e}$ is defined with  $c=0$),
$$\P(X_1^{\e,\e x_0}\in B(\e x_0+e,\alpha))>0.$$
Consequently
$$\P\Bigl(\tilde{X}_{1/\e^2}^{x_0}\in B(x_0+\frac{e}{\e},\frac{\alpha}{\e})\Bigr)>0.$$
It then follows from Stroock-Varadhan's support theorem that there exists a control
$u\in L^2(0,1/\e^2;\R^m)$ such that the corresponding $z$-trajectory, lifted to
$\R^d$, satisfies
$$\bar z(0)=x_0,\quad |\bar z(1/\e^2)-(x_0+e/\e)|\le\alpha/\e.$$
Now from \bref{h2}, there exist $t\le t_0$ and $u\in L^2(\frac{1}{\e^2},\frac{1}{\e^2}+t;\R^m)$,
such that $\bar z(\frac{1}{\e^2}+t)=x_0+g(\gamma)$, where
$$\gamma=\Bigl\{z(s),\ 0\le s\le \frac{1}{\e^2}+t\Bigr\}\in\Gamma.$$
We have constructed a loop $\gamma$ such that
$$\Bigl|g(\gamma)-\frac{e}{\e}\Bigr|\le\frac{\alpha}{\e}+ C,$$
for some universal constant $C$ (see the remark following assumption \bref{h2}). Multiplying by $\e$,
we get
$$|\e g(\gamma)-e|\le\alpha+\e C,$$
and consequently
$$\left|\frac{g(\gamma)}{|g(\gamma)|}-e\right|\le 2(\alpha+\e C),$$
from which the result follows, provided we let $\e\le \delta/4C$ and 
choose $\alpha=\delta/4$.
\end{proof}

\begin{proof}[of Corollary \ref{cor3}]
It follows from Theorem~\ref{th2} that linear combinations of elements of $G$ with positive coefficients are dense in $\Im A$. Since they form a cone, they must therefore coincide with $\Im A$.
Denoting the dimension of $\Im A$ by $k$, it follows from Theorem~\ref{th1} that, for any $g\in G$,
there exist $\alpha_1,\ldots,\alpha_k\in\R_+$ and linearly independent 
elements $g_1,\ldots,g_k\in G$ such that
$$-g=\sum_{i=1}^k\alpha_ig_i.$$
In fact, since $g$ and all $g_i$ have integer coordinates and since the inverse of a matrix
with integer coefficients has rational coefficients, the $\alpha_i$'s must be rational. 
Consequently, there exists $n\in\N$ such that $-ng\in G$, from
which it follows immediately that $-g\in G$. We have proved that $G$ is a group. 
\end{proof}

\begin{proof}[of Theorem \ref{th3}]
The first statement follows from the fact that any admissible trajectory
starting at $x_0$ stays in $S$.

Conversely, we first note that since $S$ is connected,
 to each loop $\rho\subset S$, we can associate 
a loop $\rho'\subset S$ starting from $x_0$, with $g(\rho')=g(\rho)$. So we might 
as well
assume that $\rho$ starts from $x_0$. Denote by $\bar S$ the lift of $S$ to $\R^d$.
It now follows that $x_0$ and $x_0+g(\rho)$ belong to the same connected component
$\bar S'$ of $\bar S$. It remains to show that there exists an admissible loop $\gamma$ such that $\bar\gamma$ joins $x_0$ and $x_0+g(\rho)$. In other words, it remains to show that $g(\rho) \in G$.

For that sake, first notice that
$$\bar S=\cup_{k\in \Z^d}C(k),$$
where
$C(k)$ denotes the set of accessible points from
$x_0 + k$ by the solution of the controlled ODE \eqref{control} lifted to $\R^d$. 
Consequently we have that
 $$\bar S'=\cup_{k\in \Z^d,\ x_0+k \in \bar S'}C(k).$$
Recall that $x_0+g(\rho)\in\bar S'$. Since each $C(k)$ is open, there exists $n > 0$ and
a chain $f_1,\ldots,f_n$ of points in $\Z^d$ such that, with $f_0=0$ and $f_{n+1}=g(\rho)$, 
$$C(f_i)\cap C(f_{i+1})\not=\emptyset,\quad 0\le i\le n.$$
Pick a point $y\in C(f_i)\cap C(f_{i+1})$. There exist a control such that the
solution of \eqref{control} (lifted to $\R^d$) starting from $x_0+f_i$ reaches $y$ in finite time, and another one such that the solution starting from $x_0+f_{i+1}$ also reaches $y$
in finite time. It follows from assumption \bref{h2} that there exists a control such that
the solution starting from $y$ reaches some $x_0+f$ in finite time. Hence $f-f_i$,
$f-f_{i+1}\in G$. Since $G$ is a group, $f_{i+1}-f_i\in G$, for all $0\le i\le n$.
We have proved that $g(\rho)\in G$.
\end{proof}

\section{Homogenization of an elliptic PDE}

In this section, we show how to apply the homogenization results obtained in the previous sections
to the elliptic homogenization problem \eref{ellip}.
Let $D$ be a bounded domain in $\R^d$ with a $\CC^1$ boundary, and define
\[
\tau^\eps=\inf\{t\ge0,\, X^\eps_t\not\in\overline{D}\}.
\]
Let $\alpha\ge0$ be such that for all $x\in D$,
\begin{equation}\label{unif}
\sup_{\eps>0}\E_x\exp(\alpha\tau^\eps)<\infty.
\end{equation}
We assume that $f\in\CC(\R^{2d})$, that it is periodic with respect to its second variable, and that there exists $\delta>0$
such that
\begin{equation}\label{borne}
f(x,y)\le (\alpha-\delta)^+,\ \forall x\in\R^d,\ y\in\T.
\end{equation}
\begin{remark}
We will need conditions \eqref{unif} and \eqref{borne}
in order to deduce some uniform integrability. Condition
\eqref{unif} must be checked in each particular example. 
It is always satisfied with $\alpha=0$ (unless one is in the case $A=0$), in which case
\eqref{borne} requires that $f(x,y)<0$, $x\in\R^d,\ y\in\T$. This 
condition can be relaxed only if \eqref{unif} is satisfied with some $\alpha>0$.
\end{remark}
The solution of the elliptic PDE
\begin{equation*}\left\{
\begin{aligned}
{\mathcal L}_\eps u^\eps(x)+f\left(x,\frac{x}{\eps}\right)u^\eps(x)&=0,\, x\in D,\\
u^\eps(x)&=g(x),\, x\in\partial D,
\end{aligned}\right.
\end{equation*}
is then given by the Feynman-Kac formula
\[
u^\eps(x)=\E_x\left[g(X^\eps_{\tau^\eps})\exp\left(\int_0^{\tau^\eps}
f(X^\eps_s,\frac{X^\eps_s}{\eps})\,ds\right)\right].
\]
We define as before
\begin{equation*}
\begin{aligned}
A&=\int_{\T}(I+\nabla\hat{b})a(I+\nabla\hat{b})^{*}(x)\ \mu(dx),\\
C&=\int_{\T}(I+\nabla\hat{b})c(x)\ \mu(dx),\\
D(x)&=\int_{\T}f(x,y)\ \mu(dy),\\
{\mathcal L}&=\frac{1}{2}A_{ij}\frac{\partial^2}{\partial
x_i\partial x_j}+ C_i\frac{\partial}{\partial x_i}.
\end{aligned}
\end{equation*}
From Theorem \ref{th:homo},
$ X^\eps\Rightarrow X$ as $\eps\to0$, where
$$X_t=x+Ct+A^{1/2}W_t,\ t\ge0.$$
We assume here moreover that
\begin{equation}\label{nondeg}
A \,\text{ is strictly positive definite.} \end{equation}
 Then
\[
\tau=\inf\{t\ge0,\, X_t\not\in\overline{D}\}
\]
is an a.s.  continuous function of the limiting trajectory
$\{X_t\}$. It follows from Proposition \ref{pr:conv-erg} that
\[
\int_0^t f\left(X^\eps_s,\frac{X^\eps_s}{\eps}\right)\,ds\Rightarrow \int_0^t D(X_s)\,ds.
\]
From \eqref{unif} et \eqref{borne} we deduce the necessary uniform
integrability in order to establish the
\begin{theorem}
\label{conv_ellip} Under the conditions \bref{h1}, \bref{h2},
\bref{h4}, \bref{h5}, \eqref{unif}, \eqref{borne} and \eqref{nondeg},
\[
u^\eps(x)\to\E_x\left[g(X_{\tau})\exp\left(\int_0^\tau D(X_s)\,ds
\right)\right],\; \text{as}\; \eps\to0,
\]
where $u(x):=\E_x\left[g(X_{\tau})\exp\left(\int_0^\tau D(X_s)\,ds
\right)\right]$ is the solution of the elliptic PDE
\begin{equation*}\left\{
\begin{aligned}
{\mathcal L} u(x)+D(x) u(x)&=0,\, x\in D;\\
u(x)&=g(x),\, x\in\partial D,
\end{aligned}\right.
\end{equation*}
at least in the viscosity sense.
\end{theorem}
\begin{remark}
The assumption that $A$ be nondegenerate is not really necessary for our
purpose. All we need is that for almost all (with respect to the law
of $X^x_\tau$) points $x\in\partial D$, $\langle An(x),n(x)\rangle>0$,
where $n(x)$ denotes the normal at $x$ to $\partial D$. Depending on the geometry of $D$, this can be true even if dim$(\text{Im} A)=1$
and $d>1$.
\end{remark}
\section{Homogenization of a parabolic PDE}
Assume that $a, b, c$ satisfy the above assumptions. Let $e\in \CC^1(\R^d, \R)$, 
$f \in \CC_b(\R^d\times\R^d,\R)$,
$e$ is periodic, $f$ is periodic with respect to its
second argument, and $g\in \CC(\R^d)$ grows at most
polynomially at infinity. For each $\eps>0$, we consider the PDE:
\begin{equation}\label{edp}
\left\{
\begin{aligned}
\frac{\partial u^{\eps}}{\partial t}(t, x)&={\mathcal L}_{\eps}u^{\eps}(t,x)
 + \Bigl(\frac{1}{\eps}e(\frac{x}{\eps})+f(x,\frac{x}{\eps})\Bigr)u^{\eps}(t,
x)\\  u^{\eps}(0, x)&=g(x), x\in\R^d.
\end{aligned}\right.
\end{equation}
 We assume that
 \begin{equation}\label{ezero}
\int_\T e(x)\mu(dx)=0.
\end{equation} Define
$$Y^{\eps}_{t}=\int_{0}^{t}\left[\frac{1}{\eps}e\left(\frac{X^{\eps}_{s}}{\eps}\right)
+f\left(X^{\eps}_{s},\frac{X^{\eps}_{s}}{\eps}\right)\right]\,ds.$$ Then the solution of \eqref{edp}
is given by
\begin{equ}[e:solu]
u^{\eps}(t, x)= \E[g(X^{\eps}_{t})\exp(Y^{\eps}_{t})]\;,
\end{equ}
where
$X^{\eps}_{t}$ is the solution of the SDE \eqref{sde}.

Define $$\hat{e}(x)=\int_{0}^{\infty}\E_{x}[e(\tilde{X}_{t})]\,dt,$$
$$\hat{b}_{i}(x)= \int_{0}^{\infty}\E_{x}[b_{i}(\tilde{X}_{t})]\,dt,\quad i=1,\ldots,d$$
the weak sense solutions of the Poisson equations
$$L\hat{e}(x)+e(x)=0,$$ $$L\hat{b}_{i}(x)+b_{i}(x)=0,\quad i=1,\dots, d .$$

We let \begin{equation}\label{hom_coef}
\left\{
\begin{aligned} A&=\!\int_{\T}\!\!(I+\nabla\hat{b})a(I+\nabla\hat{b})^{*}(x)\mu(dx)\\
C&=\!\int_{\T}\!\!(I+\nabla\hat{b})(c+a\nabla\hat{e})(x)\mu(dx)\\
D(x)&=\!\int_{\T}\!\!(\frac{1}{2}
\nabla\hat{e}^{*}\! a\nabla\hat{e}\!+\!f(x,\cdot)\!+\!\nabla\hat{e}c)(y)\mu(dy)
\end{aligned}
\right.
\end{equation}
Then, defining again $X^x_t=x+Ct+A^{\frac{1}{2}}W_{t}$, 
$$u(t, x)= 
\E\bigl[g(X^x_t)e^{\int_0^tD(X^x_s)\,ds}\bigr]$$ 
is the
solution (at least in the viscosity sense) of the parabolic PDE
\begin{equation*}\left\{
\begin{aligned}
\frac{\partial u}{\partial t}(t,x)&=\frac{1}{2}A_{ij}\frac{\partial^{2}u}{\partial x_{i}\partial
x_{j}}(t, x)+C_{i}\frac{\partial u}{\partial x_{i}}(t, x)+D(x) u(t,x), \\
 u(0, x)&=g(x),  \  x\in\R^d.
\end{aligned}\right.
\end{equation*}

\begin{theorem}\label{conv_edp} Under the conditions \bref{h1},
\bref{h2}, \bref{h4}, \bref{h5}, the smoothness of $e$, the
boundedness and continuity of $f$, the continuity and growth condition on $g$ and \eqref{ezero},
 for all $t\geq 0$, $x\in\R^d$,
 $$u^{\eps}(t,x)\to u(t, x)$$ as $\eps\to 0.$
\end{theorem}
For the proof of Theorem \ref{conv_edp},
we need a result which is proved exactly as Lemma \ref{le:ito}, namely
\begin{lemma}\label{le:ito2}
Define
$$ \hat{Y}^\eps_t=Y^\eps_t+\eps\hat{e}\bigl(\frac{X^\eps_t}{\eps}\bigr).$$
Then the following holds for each $\eps>0$ and $t>0$~:
$$ \hat{Y}^\eps_t=\eps\hat{e}\left(\frac{x}{\eps}\right)
+\int_0^t[f+\nabla\hat{e}c]\left(X^{\eps}_{s},\frac{X^\eps_s}{\eps}\right)\,ds
+\int_0^t\nabla\hat{e}\sigma\left(\frac{X^\eps_s}{\eps}\right)\,dW_s.$$
\end{lemma}
We can now proceed with the
\begin{proof}[of Theorem \ref{conv_edp}]
We first define a
new probability $\tilde{\P}$ by the formula
\begin{equation*}
\frac{d\tilde\P}{d\P}\Big|_t=\exp\left(\int_0^t\nabla\hat{e}\sigma\left(\frac{X^\eps_s}{\eps}\right)\,dW_s
 -\frac{1}{2}\int_0^t \scal{\nabla\hat{e}, a \nabla\hat{e}}\left(\frac{X^\eps_s}{\eps}\right)\,ds\right).
\end{equation*}
We next remark that it follows from \eref{e:solu} that the behaviour as $\eps\to0$ of $u^\eps(t,x)$ is the
same as that of 
$$\hat{u}^\eps(t,x):=\E\bigl[g(\hat{X}^{\eps}_{t})\exp\bigl(\hat{Y}^\eps_t - \eps \hat e(x/\eps)\bigr)\bigr].$$
The definition of $\tilde\P$ then yields 
$$\hat{u}^\eps(t,x)=\tilde\E\Bigl[g(\hat{X}^{\eps}_{t})\exp\int_0^t\tilde{f}
\left(X^{\eps}_{s},\frac{X^\eps_s}{\eps}\right)\,ds\Bigr]\;,$$
where we defined
$$\tilde{f}=f+\nabla\hat{e}c+\frac{1}{2} \scal{\nabla\hat{e},a \nabla\hat{e}}.$$ On the other hand, it follows from Girsanov's theorem that
\begin{equation*}
\begin{split}
\hat{X}^\eps_t &=x+\eps\hat{b}\left(\frac{x}{\eps}\right)
+\int_0^t(I+\nabla\hat{b})(c + a\nabla\hat{e})\left(\frac{X^\eps_s}{\eps}\right)\,ds\\&\  +
\int_0^t(I+\nabla\hat{b})\sigma\left(\frac{X^\eps_s}{\eps}\right)d\tilde W_s \;,
\end{split}
\end{equation*}
where $\tilde W_t=W_t-\int_0^t\sigma^t \nabla\hat{e}\left(\frac{X^\eps_s}{\eps}\right)\,ds$ is a 
Brownian motion under $\tilde\P$.

An obvious adaptation of the proof of Theorem \ref{th:homo} shows
 that 
$\hat{X}^{x,\eps}\Rightarrow X^x$, where $X^x_t=x+Ct+A^{1/2}\tilde W_t$,
$C$ is given by \eref{hom_coef}, and $\{\tilde W_t\}$ is a
Brownian motion under $\tilde \P$. Since it follows
from Proposition \ref{pr:conv-erg} that
$$\int_0^t\tilde{f}\Bigl(X^{\eps}_{s},\frac{X^\eps_s}{\eps}\Bigr)\,ds\Rightarrow \int_0^t D(X^x_s)\,ds
$$
as $\eps\to0$, the result follows.
\end{proof}

\section{Examples and counterexamples}

This section provides several examples of degenerate diffusions for which our results apply.
In each case, we furthermore give a characterisation of the range of the effective diffusion matrix $A$.

The figures should be interpreted as follows.
The little black arrows show the vector field $b$. The shaded dark grey regions denote the
points where the strong H\"ormander condition holds.

Our first example (which is the only one that we are going to work out in some detail) is
a typical example of the type of diffusions for which our conditions \bref{h1}--\bref{h5} apply.
Observe first that the diffusion $\{\tilde{X}_t,\:t\ge0\}$ on $\T$ has the measure
with density $p$ as invariant measure, and \bref{h5} is satisfied
if and only if there exists a mapping $H$  from $\T$ into the set of
$d\times d$ antisymmetric matrices, such that
\[
b_i(x)=\frac{1}{2p(x)}\sum_j\frac{\partial}{\partial x_j}\left(
p a_{ij}+H_{ij}\right)(x),\:1\le i\le d,\:x\in\T.
\]

Let $\alpha:\T
\to[0,1]$ be a smooth function such that the set $\CA = \{x \,|\, \alpha(x)=0\}$ has a 
finite 
number of bounded connected components.
For $\rho > 0$, we define the set $\CA_\rho:=\{x,\: d(x,\CA)\le\rho\}$
and we assume that there exists $\rho > 0$ such that $\CA_\rho$ intersects neither $\{x_1 = 0\}$
nor $\{x_1 = 1/2\}$.

We choose our diffusion matrix $\sigma$ to be given by $\sigma(x)= \alpha(x) I$,
and we let
\[
H=\begin{pmatrix}
0&0&\cdots&0&(2\pi)^{-1}\cos(2\pi x_1)\\
0&0&\cdots&0&0\\
\vdots&\vdots&\cdots&\vdots&\vdots\\
0&0&\cdots&0&0\\
-(2\pi)^{-1}\cos(2\pi x_1)&0&\cdots&0&0
\end{pmatrix}\;.
\]
Hence, if we choose 
\[
2b(x)=\nabla(\alpha^2)(x)+\begin{pmatrix}0\\ \vdots\\0\\ \sin(2\pi x_1)
\end{pmatrix}\;,
\]
then the Lebesgue measure is invariant for the process $\tilde X^x_t$.
The following picture shows the vector field given by $b$, together with an example
of a set $\CA$ (the white disks) that satisfies our conditions.
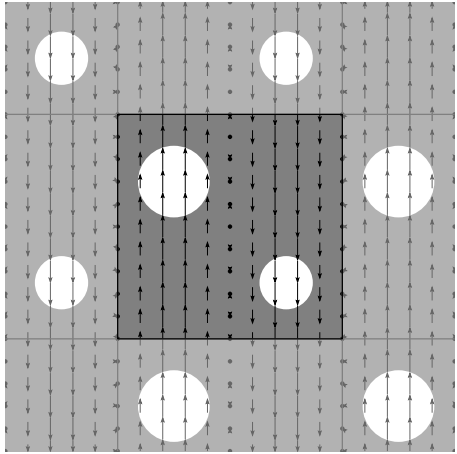
\begin{equ}
\mhpastefig{etienne}
\end{equ}
Assumptions \bref{h1} and \bref{h5} are satisfied by construction. Assumption \bref{h2} is satisfied since,
by taking the trivial control $u=0$, there exists a time $t$ such that every solution to the controlled
system reaches the complement of $\CA$ before time $t$. The same argument shows that \bref{h4}
is satisfied as well.

Our second example is one where the diffusion matrix of the limiting Brownian motion
is non-degenerate even though the area in which the strong H\"ormander condition holds
does not intersect the boundaries of the fundamental domain. In particular, the argument of \cite{JKO}
mentioned in Remark~\ref{rem:JKO} does not cover this situation:
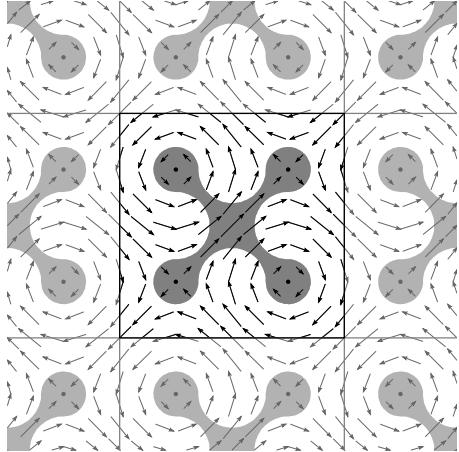
\begin{equ}
\mhpastefig{zero}
\end{equ}
In this case, it is easy to see that every point of the torus can be reached by the diffusion on the torus, 
so that $\mu$ has full support and therefore $A$ has full rank.

Our third example is in a way the opposite of the first one. It shows that it is possible
for the limiting Brownian motion to have zero diffusion coefficient, even though the area
in which the strong H\"ormander condition holds stretches over the whole space
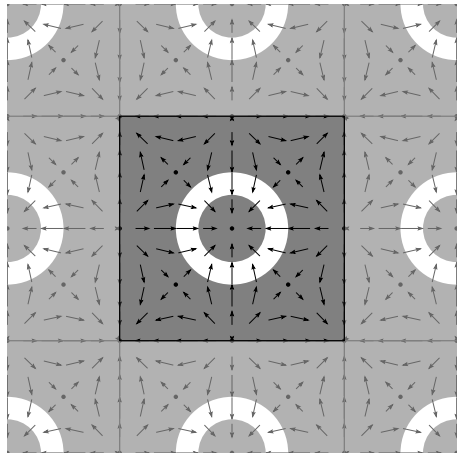
\begin{equ}
\mhpastefig{degenerate}
\end{equ}
In this case, the support $S$ of the invariant measure is given by the gray disk, so that 
$g(\gamma) = 0$ for every loop in $S$.

We finally give an example where the effective diffusivity degenerates in one direction.
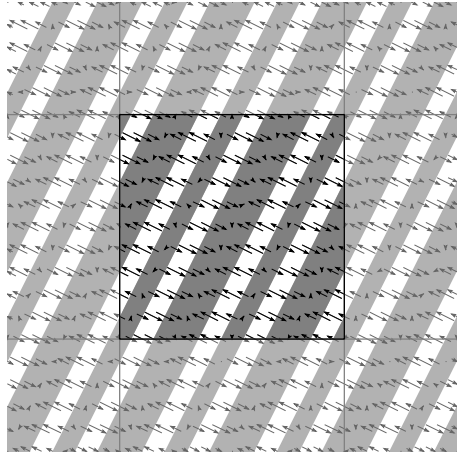
\begin{equ}
\mhpastefig{skewed}
\end{equ}
In this particular case, the range of $A$ is the one-dimensional subspace of $\R^2$
spanned by the vector $(1,2)$.

\bibliography{refs}
\bibliographystyle{Martin}

\end{document}